\input amstex
\input amsppt.sty
\magnification=\magstep1
\hsize=30truecc
\vsize=22.2truecm
\baselineskip=16truept
\NoBlackBoxes
\TagsOnRight \pageno=1 \nologo
\def\Z{\Bbb Z}
\def\N{\Bbb N}

\def\l{\left}
\def\r{\right}
\def\bg{\bigg}
\def\({\bg(}
\def\[{\bg\lfloor}
\def\){\bg)}
\def\]{\bg\rfloor}
\def\t{\text}
\def\f{\frac}

\def\p{\ (\roman{mod}\ p)}

\def\bi{\binom}
\def\eq{\equiv}

\def\ls{\leqslant}
\def\gs{\geqslant}
\def\mo{\roman{mod}}

\def\da{\delta}

\def\Proof{\noindent{\it Proof}}

\def\Remark{\medskip\noindent{\it  Remark}}

\def\Ack{\medskip\noindent {\bf Acknowledgments}}
\hbox {Sci. China Math. 54(2011), no.\,12, 2509--2535.}
\bigskip
\topmatter
\title Super congruences and Euler numbers\endtitle
\author Zhi-Wei Sun\endauthor
\leftheadtext{Zhi-Wei Sun} \rightheadtext{Super congruences and
Euler numbers}
\affil Department of Mathematics, Nanjing University\\
 Nanjing 210093, People's Republic of China
  \\  zwsun\@nju.edu.cn
  \\ {\tt http://math.nju.edu.cn/$\sim$zwsun}
\endaffil
\abstract Let $p>3$ be a prime. A $p$-adic congruence is called a super congruence if it happens to hold modulo some higher power of $p$.
The topic of super congruences is related to many fields including Gauss and Jacobi sums and hypergeometric series.
We prove that
$$\align\sum_{k=0}^{p-1}\f{\bi{2k}k}{2^k}\eq& (-1)^{(p-1)/2}-p^2E_{p-3}\ (\mo\ p^3),
\\\sum_{k=1}^{(p-1)/2}\f{\bi{2k}k}k\eq&(-1)^{(p+1)/2}\f 83 pE_{p-3}\ (\mo\ p^2),
\\\sum_{k=0}^{(p-1)/2}\f{\bi{2k}k^2}{16^k}\eq&(-1)^{(p-1)/2}+p^2E_{p-3}\ (\mo\ p^3),
\endalign$$
where $E_0,E_1,E_2,\ldots$ are Euler numbers. Our new approach is of combinatorial nature.
We also formulate many conjectures concerning super congruences and relate most of them to Euler numbers or Bernoulli numbers.
Motivated by our investigation of super congruences, we also raise a conjecture on 7 new series
for $\pi^2$, $\pi^{-2}$ and the constant $K:=\sum_{k=1}^\infty(\f k3)/k^2$ (with $(-)$ the Jacobi symbol), two of which are
$$\sum_{k=1}^\infty\f{(10k-3)8^k}{k^3\bi{2k}k^2\bi{3k}k}=\f{\pi^2}2
\ \ \ \t{and}\ \ \ \sum_{k=1}^\infty\f{(15k-4)(-27)^{k-1}}{k^3\bi{2k}k^2\bi{3k}k}=K.$$
\endabstract
\thanks 2010 {\it Mathematics Subject Classification}.\,Primary 11B65;
Secondary 05A10, 05A19, 11A07, 11B68, 11E25, 11F20, 11M06, 11S99, 33C20.
\newline\indent The initial version of this paper was posted to arXiv as a preprint in Jan. 2010 with the code {\tt arXiv:1001.4453}.
\newline\indent Supported by the National Natural Science
Foundation (grant 10871087) and the Overseas Cooperation Fund (grant 10928101) of China.
\endthanks
\endtopmatter
\document

\heading{1. Introduction}\endheading

Let $p$ be an odd prime. Clearly
$$\bi{2k}k=\f{(2k)!}{(k!)^2}\eq0\ (\mo\ p)\quad\t{for every}\ k=\f{p+1}2,\ldots,p-1.$$
After a series of work on combinatorial congruences involving central binomial coefficients (cf. [PS], [ST1] and [ST2]),
Z. W. Sun [S10b] determined
$\sum_{k=0}^{p-1}\bi{2k}k/m^k$ modulo $p^2$ for any integer $m\not\eq0\ (\mo\ p)$. In particular, he showed that
$$\sum_{k=0}^{p-1}\f{\bi{2k}k}{2^k}\eq(-1)^{(p-1)/2}\ (\mo\ p^2)$$
and conjectured that there are no odd composite numbers $n>1$
satisfying the congruence
$\sum_{k=0}^{n-1}\bi{2k}k/2^k\eq(-1)^{(n-1)/2}\ (\mo\ n^2)$.
He also searched those {\it exceptional primes} $p$ such that
$\sum_{k=0}^{p-1}\bi{2k}k/2^k\eq (-1)^{(p-1)/2}\ (\mo\ p^3)$ and only found
two such primes: $149$ and $241$.

Let $p$ be an odd prime.
A $p$-adic congruence is said to be a {\it super congruence} if it happens to hold modulo some higher power of $p$.
In 2003 Rodriguez-Villegas [RV] conjectured  22 super congruences
via his analysis of the $p$-adic analogues of Gaussian hypergeometric series and
the Calabi-Yau manifolds.
The most elegant one of the 22 super congruences is as follows:
$$\sum_{k=0}^{p-1}\f{\bi{2k}k^2}{16^k}\eq(-1)^{(p-1)/2}\ (\mo\ p^2),$$
which was later proved by E. Mortenson [M03a] via the $p$-adic $\Gamma$-function and the Gross-Koblitz formula.
See also K. Ono's book [O] and the papers [AO], [K], [MO], [M03b], [M05], [M08], [Mc] and [OS]
for such advanced approach to super congruences.
Recently the author's twin brother Z. H. Sun, as well as R. Tauraso [T], gave a simple proof of the last
super congruence, and Z. H. Sun also proved the author's conjectural congruence
$$\sum_{k=0}^{p-1}\f{k\bi{2k}k^2}{16^k}\eq\f{(-1)^{(p+1)/2}}4\ (\mo\ p^2)$$
via a combinatorial identity. Note that Stirling's formula
$n!\sim \sqrt{2\pi n}(n/e)^n$ implies that
$$\lim_{k\to+\infty}\f{k\bi{2k}k^2}{16^k}=\f1{\pi}.$$

Surprisingly, we find that the above topics  are related to Euler numbers.

Recall that Euler numbers $E_n\ (n\in\N=\{0,1,2,\ldots\})$ are integers defined by
$$E_0=1,\ \t{and}\ \sum^n\Sb k=0\\2\mid k\endSb \bi nk E_{n-k}=0\ \ \ \t{for}\ n\in\Z^+=\{1,2,3,\ldots\}.$$
It is well known that $E_{2n+1}=0$ for all $n\in\N$ and
$$\sec x=\sum_{n=0}^\infty(-1)^n E_{2n}\f{x^{2n}}{(2n)!}\ \ \t{for}\ |x|<\f{\pi}2.$$

Now we state the main results and key conjectures in this paper.

\proclaim{Theorem 1.1} Let $p$ be an odd prime. Then
$$\sum_{k=0}^{p-1}\f{\bi{2k}k}{2^k}\eq (-1)^{(p-1)/2}-p^2E_{p-3}\ (\mo\ p^3).\tag1.1$$
If $p>3$, then
$$\sum_{k=1}^{(p-1)/2}\f{\bi{2k}k}k\eq (-1)^{(p+1)/2}\,\f 83pE_{p-3}\ (\mo\ p^2)\tag1.2$$
and
$$\sum_{k=1}^{(p-1)/2}\f1{k^2\bi{2k}k}\eq(-1)^{(p-1)/2}\,\f 43E_{p-3}\ (\mo\ p).\tag1.3$$
We also have
$$\sum_{k=1}^{(p-1)/2}\f{4^k}{k^2\bi{2k}k}\eq(-1)^{(p-1)/2}\,4E_{p-3}\ (\mo\ p)\tag1.4$$
and
$$\sum_{p/2<k<p}\f{\bi{2k}k}{k4^k}\eq (-1)^{(p-1)/2}\,2pE_{p-3}\ (\mo\ p^2).\tag1.5$$
\endproclaim

\Remark\ 1.1. By (1.1) those exceptional primes are just those odd primes $p$ with $p\mid E_{p-3}$;
the two exceptional primes $149$ and $241$ offer the main clue to our discovery of (1.2)-(1.5).
Also, Sun and Tauraso [ST1] showed that
$$\sum_{k=1}^{p-1}\f{\bi{2k}k}k\eq\f89 p^2B_{p-3}\ (\mo\ p^3)$$
for any prime $p>3$, where
$B_0,B_1,B_2,\ldots$ are Bernoulli numbers. It is remarkable that
$$\sum_{k=1}^\infty \f1{k^2\bi{2k}k}=\f{\pi^2}{18}\ \ \t{and}\ \  \sum_{k=1}^\infty \f{4^k}{k^2\bi{2k}k}=\f{\pi^2}2$$
(see [Po, (3)] and problem 44b of [St, Chapter 1] for the first series, and [Ma] and [Sp] for the second series),
which were even known in the nineteenth century.
Tauraso ([T1],[T2]) showed that
$$\sum_{k=1}^{p-1}\f{\bi{2k}k}{k4^k}\eq-\sum_{k=1}^{(p-1)/2}\f1k\ (\mo\ p^3)\ \
\t{and}\ \
\sum_{k=1}^{p-1}\f1{k^2\bi{2k}k}\eq \f1{3p}\sum_{k=1}^{p-1}\f1k\ (\mo\ p^3)$$
 for any prime $p>5$.

\medskip

Recall that harmonic numbers are those integers
$$H_n:=\sum_{0<k\ls n}\f 1k\ \ \ (n=0,1,2,\ldots).$$
It is known (cf. [S1]) that
$$\f{H_{p-1}}{p^2}\eq-\f{B_{p-3}}3\ (\mo\ p)\ \ \t{and}\ \ \f5{p^2}\sum_{k=1}^{p-1}\f1{k^3}\eq-6B_{p-5}\ (\mo\ p)$$
for any prime $p>3$.

Now we present our first conjecture.

\proclaim{Conjecture 1.1} Let $p>3$ be a prime. Then
$$\sum_{k=1}^{p-1}\f{4^k}{k^2\bi{2k}k}+\f{4q_p(2)}p\eq-2q^2_p(2)+pB_{p-3}\ (\mo\ p^2)$$
and $$p\sum_{k=1}^{p-1}\f{2^k}{k\bi{2k}k}\eq \l(\f{-1}p\r)-1-p\,q_p(2)+p^2E_{p-3}\ (\mo\ p^3),$$
where $(-)$ denotes the Jacobi symbol and $q_p(2)$ stands for the Fermat quotient $(2^{p-1}-1)/p$.
Also, $$\sum_{k=1}^{p-1}\f{\bi{2k}k}{k^3}\eq\f 23B_{p-3}\ (\mo\ p),$$
and furthermore
$$\sum_{k=1}^{p-1}\f{\bi{2k}k}{k^3}\eq-\f 2{p^2}H_{p-1}-\f{13}{27}\sum_{k=1}^{p-1}\f1{k^3}\ (\mo\ p^4)\ \ \t{if}\ p>7.$$
When $p>5$ we have
$$\align\sum_{k=1}^{(p-1)/2}\f{(-1)^k}{k^3\bi{2k}k}\eq&-2B_{p-3}\ (\mo\ p),
\\\sum_{k=1}^{(p-1)/2}\f{(-1)^k}{k^2}\bi{2k}k\eq&\f{56}{15}pB_{p-3}\ (\mo\ p^2),
\\\sum_{k=1}^{p-1}\f1{k^4\bi{2k}k}-\f{H_{p-1}}{p^3}\eq&-\f 7{45}pB_{p-5}\ (\mo\ p^2).
\endalign$$
\endproclaim
\Remark\ 1.2. It is known that
$$\sum_{k=1}^\infty\f{2^k}{k\bi{2k}k}=\f{\pi}2,\ \ \sum_{k=1}^\infty\f{(-1)^k}{k^3\bi{2k}k}=-\f25\zeta(3)
\ \ \t{and}\ \ \sum_{k=1}^\infty\f1{k^4\bi{2k}k}=\f{17}{36}\zeta(4).$$
Tauraso [T2] determined $\sum_{k=1}^{p-1}(-1)^k/(k^3\bi{2k}k)$ and $\sum_{k=1}^{p-1}(-1)^k\bi{2k}k/k^2$
modulo $p^2$ (for any prime $p>5$) in terms of $H_{p-1}$.

\proclaim{Theorem 1.2} Let $p>3$ be a prime.
Then
$$\align\sum_{k=0}^{(p-1)/2}\f{\bi{2k}k}{8^k}\eq&\l(\f 2p\r)+\l(\f{-2}p\r)\f{p^2}4E_{p-3}\ (\mo\ p^3);\tag1.6$$
\\\sum_{k=0}^{(p-1)/2}\f{\bi{2k}k^2}{16^k}\eq&(-1)^{(p-1)/2}+p^2E_{p-3}\ (\mo\ p^3),\tag1.7$$
\\\sum_{k=0}^{(p-1)/2}\f{k\bi{2k}k^2}{16^k}\eq&\f{(-1)^{(p+1)/2}}4+\f{p^2}4(1-E_{p-3})\ (\mo\ p^3);\tag1.8$$
\\\sum_{p/2<k<p}\f{\bi{2k}k^2}{16^k}\eq&-2p^2E_{p-3}\ (\mo\ p^3),\tag1.9
\\\sum_{p/2<k<p}\f{k\bi{2k}k^2}{16^k}\eq&\f{p^2}2 E_{p-3}\ (\mo\ p^3).\tag1.10
\endalign$$
Furthermore,
$$(-1)^{(p-1)/2}\sum_{k=0}^{(p-1)/2}\f{\bi{2k}k^2}{16^k}\eq 1-\f38p\sum_{k=1}^{(p-1)/2}\f{\bi{2k}k}k\ (\mo\ p^4)\tag1.11$$
and
$$\aligned\sum_{k=0}^{(p-1)/2}\f{k\bi{2k}k^2}{16^k}\eq&\f{(-1)^{(p+1)/2}}4+\f{p^2}4(2^{p}-1)
\\&+(-1)^{(p-1)/2}\f 3{32}p\sum_{k=1}^{(p-1)/2}\f{\bi{2k}k}k\ (\mo\ p^4).
\endaligned\tag1.12$$
\endproclaim

\Remark\ 1.3.  The reason why we don't include (1.6) in Theorem 1.1
is that its proof is similar to that of (1.11) and (1.12).
For any prime $p>3$,  R. Osburn and C. Schneider [OS] used Jacobi sums and the $p$-adic $\Gamma$-function
to prove that
$$(-1)^{(p-1)/2}\sum_{k=0}^{(p-1)/2}\f{\bi{2k}k^2}{16^k}\eq 1-\f38p\sum_{k=1}^{(p-1)/2}\f{\bi{2k}k}k\ (\mo\ p^3).$$
The author observed that a combination of (1.11) and (1.12) yields that
$$\sum_{k=0}^{(p-1)/2}\f{4k+1}{16^k}\bi{2k}k^2\eq p^2(2^p-1)\ (\mo\ p^4)$$
for any prime $p>3$. After reading this, Tauraso noted the following identity
$$\sum_{k=0}^n\f{4k+1}{16^k}\bi{2k}k^2=\f{(n+1)^2}{16^n}\bi{2n+1}n^2=\f{(2n+1)^2}{16^n}\bi{2n}n^2,$$
which can be easily proved by induction. This identity
implies the author's following observation:
$$\sum_{p/2<k<p}\f{4k+1}{16^k}\bi{2k}k^2\eq 6p^2(1-2^{p-1})\ (\mo\ p^4)$$
for each odd prime $p$.

\proclaim{Conjecture 1.2} Let $p$ be an odd prime and let $a\in\Z^+$.
If $p\eq1\ (\mo\ 4)$ or $a>1$, then
$$\sum_{k=0}^{\lfloor \f34p^a\rfloor}\f{\bi{2k}k}{(-4)^k}\eq\l(\f2{p^a}\r)\ (\mo\ p^2)
\ \ \t{and}\ \ \sum_{k=0}^{\lfloor
\f34p^a\rfloor}\f{\bi{2k}k^2}{16^k}\eq \l(\f{-1}{p^a}\r)\ (\mo\
p^3).\tag1.13$$ If $p>3$, and $p\eq1,3\ (\mo\ 8)$ or $a>1$, then
$$\sum_{k=0}^{\lfloor \f r8p^a\rfloor}\f{\bi{2k}k^2}{16^k}\eq \l(\f{-1}{p^a}\r)\ (\mo\ p^3)
\ \ \ \t{for}\ r=5,7.\tag1.14$$
\endproclaim

Here is our third theorem.

\proclaim{Theorem 1.3} Let $p$ be a prime and let $a\in\Z^+$. Then
$$\f1{p^a}\sum_{k=0}^{p^a-1}(21k+8)\bi{2k}k^3\eq 8+16p^{3}B_{p-3} \ (\mo\ p^{4}),\tag1.15$$
where $B_{-1}$ is regarded as zero.
\endproclaim

\Remark\ 1.4. In [S11a] the author conjectured that for any odd prime $p$ we have
$$\sum_{k=0}^{p-1}\bi{2k}k^3\eq \cases4x^2-2p\ (\mo\ p^2)&\t{if}\ (\f p7)=1\ \&\ p=x^2+7y^2\ (x,y\in\Z),
\\0\ (\mo\ p^2)&\t{if}\ (\f p7)=-1,\ \t{i.e.},\ p\eq3,5,6\ (\mo\ 7).\endcases$$
Let $p$ be a prime with $(\f p7)=-1$. We also conjecture that $\f1n\sum_{k=0}^{n-1}\bi{2k}k^3$ is a $p$-adic integer
for any $n\in\Z^+$, and that
$$\sum_{k=0}^{p^a-1}\bi{2k}k^3\eq\cases0\pmod{p^{a+1}}&\t{if}\ a\in\{1,3,5,\ldots\},
\\p^a\pmod{p^{a+3-\delta_{p,3}}}&\t{if}\ a\in\{2,4,6,\ldots\},\endcases$$
where the Kronecker symbol $\da_{p,3}$ takes 1 or 0 according as $p=3$ or not.
\medskip

In a previous version of this paper, the author conjectured that
for any positive integer $n$ the arithmetic mean
$$s_n:=\f1n\sum_{k=0}^{n-1}(21k+8)\bi{2k}k^3\tag1.16$$
is always an integer divisible by $4\bi{2n}n$, and observed the recursion
$$n^3(n+1)s_{n+1}=n^4s_n+8(2n-1)^3(21n+8)\bi{2(n-1)}{n-1}^3\ \ (n=1,2,3,\ldots).$$
On Feb. 11, 2010, Kasper Andersen noted that this recurrence relation yields the following recursion for $t_n:=s_n/(4\bi{2n}n)$:
$$(4n+2)t_{n+1}-nt_n=(21n+8)\bi{2n-1}n^2\ (n=1,2,3,\ldots).$$
Then Andersen used Zeilberger's algorithm (cf. [PWZ]) to find that
$$r_n:=\sum_{k=0}^{n-1}\bi{n+k-1}k^2\ (n=1,2,3,\ldots)\tag1.17$$
 satisfies the same recursion
and hence he obtained that $t_n=r_n\in\Z$ since $t_1=r_1$.
Thanks to Andersen's discovery, we are now able to prove Theorem 1.3
which was an earlier conjecture of the author.

We guess that any integer $n>1$ satisfying $s_n\eq8\ (\mo\ n^3)$
must be a prime; this has been verified for $n\ls 10^4$.
It seems that $t_n\not\eq3\ (\mo\ 4)$, and $t_n$ is composite for all $n=3,4,\ldots$.
It is interesting to compare $s_n$ and $t_n$ with Ap\'ery numbers (cf. [Po]).

\proclaim{Conjecture 1.3}
If $p$ is a prime and $a$ is a positive integer with $p^a\eq1\ (\mo\ 3)$,
then
$$\sum_{k=0}^{\lfloor\f23p^a\rfloor}(21k+8)\bi{2k}k^3\eq8p^a\ (\mo\ p^{a+5+(-1)^p}).\tag1.18$$
Also, for each prime $p>5$ we have
$$\sum_{k=1}^{p-1}\f{21k-8}{k^3\bi{2k}k^3}+\f{p-1}{p^3}
\eq\f{H_{p-1}}{p^2}(15p-6)+\f{12}5p^2B_{p-5}\ (\mo\ p^3).\tag1.19$$
\endproclaim

It is interesting to compare Theorem 1.3 and Conjecture 1.3 with the following elegant identity
$$\sum_{k=1}^\infty\f{21k-8}{k^3\bi{2k}k^3}=\zeta(2)=\f{\pi^2}6$$
obtained by D. Zeilberger [Z] via the WZ method. In the same spirit,
we formulate the following conjecture inspired by
our observations of some congruences (see Conjectures 5.3-5.6, Remark 5.2 and Conj. 5.15(i) in Section 5).

\proclaim{Conjecture 1.4} We have
$$\sum_{k=1}^\infty\f{(10k-3)8^k}{k^3\bi{2k}k^2\bi{3k}k}=\f{\pi^2}2,\
\sum_{k=1}^\infty\f{(11k-3)64^k}{k^3\bi{2k}k^2\bi{3k}k}=8\pi^2,
\ \sum_{k=1}^\infty\f{(35k-8)81^k}{k^3\bi{2k}k^2\bi{4k}{2k}}=12\pi^2.\tag{1.20}$$
Also,
$$\sum_{k=1}^\infty\f{(15k-4)(-27)^k}{k^3\bi{2k}k^2\bi{3k}k}=-27K
\ \ \t{and}\ \  \sum_{k=1}^\infty\f{(5k-1)(-144)^k}{k^3\bi{2k}k^2\bi{4k}{2k}}=-\f{45}2K,\tag{1.21}$$
where
$$K:=L\l(2,\l(\f {\cdot}3\r)\r)=\sum_{k=1}^\infty\f{(\f k3)}{k^2}=0.781302412896486296867187429624\ldots.$$
Moreover,
$$\sum_{n=0}^\infty\f{18n^2+7n+1}{(-128)^n}\bi{2n}n^2\sum_{k=0}^n\bi{-1/4}k^2\bi{-3/4}{n-k}^2=\f{4\sqrt2}{\pi^2}\tag1.22$$
and
$$\sum_{n=0}^\infty\f{40n^2+26n+5}{(-256)^n}\bi{2n}n^2\sum_{k=0}^n\bi{n}k^2\bi{2k}k\bi{2(n-k)}{n-k}
=\f{24}{\pi^2}.\tag1.23$$
\endproclaim

One can easily check the identities in Conjecture 1.4 numerically.
Let us take the first identity in (1.20) as an example. The series converges rapidly since
$$\bi{2k}k^2\bi{3k}k\sim\f{\sqrt3}2\cdot\f{108^k}{(k\pi)^{1.5}}\quad (k\to+\infty)$$
by Stirling's formula. Via {\tt Mathematica} we find that
$$\bigg|\f2{\pi^2}\sum_{k=1}^{200}\f{(10k-3)8^k}{k^3\bi{2k}k^2\bi{3k}k}-1\bigg|<\f1{10^{227}}.$$
This provides a powerful evidence to support the first identity in (1.20).
\medskip

We will show Theorems 1.1--1.3 in Sections 2-4 respectively;
our new approach to super congruences is of combinatorial nature.
In Section 5 we will raise many new conjectures for further research.

\heading{2. Proof of Theorem 1.1}\endheading

\medskip
\noindent{\it Proof of (1.1)}. By [ST1, (2.1)], we have
$$\sum_{k=0}^{p-1}\bi{2k}k2^{p-1-k}=\sum_{k=0}^{p-1}\bi{2p}ku_{p-k},$$
where $u_0=0,\ u_1=1$ and $u_{n+1}=-u_{n-1}$ for $n=1,2,3,\ldots$. Clearly $u_{2n}=0$ and $u_{2n+1}=(-1)^n$
for all $n\in\N$. Thus
$$\sum_{k=0}^{p-1}\bi{2k}k2^{p-1-k}=\sum_{k=0}^{(p-1)/2}\bi{2p}{2k}(-1)^{(p-2k-1)/2}.\tag2.1$$
For $k=1,\ldots,(p-1)/2$, we have
$$\align \bi{2p}{2k}=&\f{2p}{2k}\bi{2p-1}{2k-1}=\f pk\prod_{j=1}^{2k-1}\f{2p-j}j=-\f pk\prod_{j=1}^{2k-1}\l(1-\f{2p}j\r)
\\\eq&-\f pk(1-2pH_{2k-1})=\f pk(1-2(1-pH_{2k-1}))
\\\eq&\f pk\(1+2\bi{p-1}{2k-1}\)=4\bi p{2k}+\f pk\ (\mo\ p^3).
\endalign$$
Thus
$$\align&(-1)^{(p-1)/2}\sum_{k=0}^{p-1}\bi{2k}k2^{p-1-k}-1
\\\eq&\sum_{k=1}^{(p-1)/2}(-1)^k\l(4\bi p{2k}+\f pk\r)
=p\sum_{k=1}^{(p-1)/2}\f{(-1)^k}k+4\sum^{p-1}\Sb k=1\\2\mid k\endSb\bi pk(-1)^{k/2}\ (\mo\ p^3).
\endalign$$

By Lehmer [L],
$$H_{(p-1)/2}\eq-2q_p(2)+p\,q^2_p(2)\ (\mo\ p^2).\tag2.2$$
In view of [S2, Corollary 3.3] we also have
$$H_{\lfloor p/4\rfloor}\eq-3q_p(2)+\f 32 p\,q^2_p(2)-(-1)^{(p-1)/2}pE_{p-3}\ (\mo\ p^2).\tag2.3$$
Therefore
$$\align \sum_{k=1}^{(p-1)/2}\f{(-1)^k}k=&\sum_{k=1}^{(p-1)/2}\f{1+(-1)^k}k-\sum_{k=1}^{(p-1)/2}\f1k
=H_{\lfloor p/4\rfloor}-H_{(p-1)/2}
\\\eq&-q_p(2)+\f p2 q_p^2(2)+(-1)^{(p+1)/2}pE_{p-3}\ (\mo\ p^2).
\endalign$$

Note also that
$$\sum^{p-1}\Sb k=1\\2\mid k\endSb\bi pk (-1)^{k/2}=\l(\f 2p\r)2^{(p-1)/2}-1$$
by [S02, (3.2)]. Combining the above we obtain
$$\align &\(\l(\f{-1}p\r)\sum_{k=0}^{p-1}\f{\bi{2k}k}{2^k}-1\)2^{p-1}+2^{p-1}-1
\\\eq&4\l(\l(\f2p\r)2^{(p-1)/2}-1\r)-p\,q_p(2)+\f{p^2}2q_p^2(2)+(-1)^{(p+1)/2}p^2E_{p-3}\ (\mo\ p^3).
\endalign$$
Observe that
$$\align &2\l(\l(\f2p\r)2^{(p-1)/2}-1\r)-p\,q_p(2)
\\=&2\l(\f 2p\r)2^{(p-1)/2}-2^{p-1}-1=-\l(\l(\f 2p\r)2^{(p-1)/2}-1\r)^2.
\endalign$$
Therefore
$$\align&\f{(\f{-1}p)\sum_{k=0}^{p-1}\bi{2k}k/2^k-1}{p^2}
\\\eq&-2\l(\f{(\f2p)2^{(p-1)/2}-1}p\r)^2+\f{q_p^2(2)}2+(-1)^{(p+1)/2}E_{p-3}\ (\mo\ p).
\endalign$$
Since
$$q_p(2)=\f{(\f2p)2^{(p-1)/2}-1}p\l(\l(\f 2p\r)2^{(p-1)/2}+1\r)\eq2\times\f{(\f2p)2^{(p-1)/2}-1}p\ (\mo\ p),$$
we finally obtain (1.1). \qed

\proclaim{Lemma 2.1} Let $p$ be an odd prime. Then, for any $k=1,\ldots,p-1$ we have
$$k\bi{2k}k\bi{2(p-k)}{p-k}\eq(-1)^{\lfloor 2k/p\rfloor-1}2p\ (\mo\ p^2).\tag2.4$$
\endproclaim
\Proof. For $k=1,\ldots,(p-1)/2$, if
$$(p-k)\bi{2(p-k)}{p-k}\bi{2k}k\eq(-1)^{\lfloor 2(p-k)/p\rfloor-1}2p=2p\ (\mo\ p^2)$$
then $$k\bi{2k}k\bi{2(p-k)}{p-k}\eq-2p=(-1)^{\lfloor2k/p\rfloor-1}2p\ (\mo\ p^2)$$
since $\bi{2(p-k)}{p-k}\eq0\ (\mo\ p)$. So it suffices to show (2.4) for any $k=(p+1)/2,\ldots,p-1$.

 Let $k\in\{(p+1)/2,\ldots,p-1\}$. Then
 $$\align \f1p\bi{2k}k=&\f1p\times\f{(2k)!}{(k!)^2}=\f1p\times\f{p!(p+1)\cdots(p+(2k-p))}{((p-1)!/\prod_{j=1}^{p-1-k}(p-j))^2}
 \\=&\f1{(p-1)!}\prod_{i=1}^{2k-p}(p+i)\times\prod_{j=1}^{p-1-k}(p-j)^2
 \\\eq&\f{(2k-p)!}{(p-1)!}((p-1-k)!)^2=\f{((p-1-k)!)^2}{\prod_{j=1}^{2(p-k)-1}(p-j)}
 \\\eq&-\f{((p-1-k)!)^2}{(2(p-k)-1)!}=-\f2{(p-k)\bi{2(p-k)}{p-k}}\eq\f2{k\bi{2(p-k)}{p-k}}\ (\mo\ p)
 \endalign$$
and hence
$$k\bi{2k}k\bi{2(p-k)}{p-k}\eq 2p\ (\mo\ p^2)$$
as desired. \qed

\Remark\ 2.1.  [T2] contains certain technique similar to Lemma 2.1.
\medskip

\proclaim{Lemma 2.2} For any $n\in\Z^+$ we have
$$\sum_{k=1}^n\f{\bi{2k}k}k=\f{n+1}3\bi{2n+1}{n}\sum_{k=1}^n\f1{k^2\bi nk^2}\tag2.5$$
and
$$\sum_{k=1}^n\f{(-1)^k}{k^2\bi nk\bi{n+k}k}=(-1)^{n-1}\(3\sum_{k=1}^n\f1{k^2\bi{2k}k}+2\sum_{k=1}^n\f{(-1)^k}{k^2}\).
\tag2.6$$
\endproclaim
\Remark\ 2.2. These two identities are known results. (2.5) is due to T. B. Staver [Sta] (see also (5.2) of [Go, p.\,50]),
and (2.6) was discovered by Ap\'ery (see [Ap] and [Po]) during his study of the irrationality of
$\zeta(3)=\sum_{n=1}^\infty1/n^3$.

\proclaim{Lemma 2.3} We have the new combinatorial identity
$$\sum_{k=1}^n\f1{k^2\bi{n+k}k}=3\sum_{k=1}^n\f1{k^2\bi{2k}k}-\sum_{k=1}^n\f1{k^2}.\tag2.7$$
\endproclaim
\Proof. Observe that
$$\align&\sum_{k=1}^n\f1{k^2\bi{n+k}k}-\sum_{k=1}^n\f1{k^2\bi{n+1+k}k}
\\=&\sum_{k=1}^n\f{\bi{n+k+1}k-\bi{n+k}k}{k^2\bi{n+k}k\bi{n+k+1}k}=\sum_{k=1}^n\f{\bi{n+k}{k-1}}{k^2\bi{n+k}k\bi{n+k+1}k}
\\=&\sum_{k=1}^n\f{n!(k-1)!}{(n+k+1)!}=\sum_{j=0}^{n-1}\f{n!j!}{(n+2+j)!}
=\f1{(n+1)(n+2)}\sum_{k=0}^{n-1}\f1{\bi{n+2+k}{n+2}}.
\endalign$$
By (2.26) of [Go, p.\,21],
$$\sum_{k=0}^m\f1{\bi{x+k}l}=\f l{l-1}\(\f1{\bi{x-1}{l-1}}-\f1{\bi{x+m}{l-1}}\)$$
for any $l\in\Z^+$. So we have
$$\sum_{k=1}^n\f1{k^2\bi{n+k}k}-\sum_{k=1}^n\f1{k^2\bi{n+1+k}k}=\f1{(n+1)(n+2)}\times\f{n+2}{n+1}\(1-\f1{\bi{(n+2)+n-1}{n+1}}\)$$
and hence
$$\align \sum_{k=1}^{n+1}\f1{k^2\bi{n+1+k}k}-\sum_{k=1}^n\f1{k^2\bi{n+k}k}
=&\f1{(n+1)^2\bi{2n+2}{n+1}}-\f1{(n+1)^2}\(1-\f1{\bi{2n+1}n}\)
\\=&\f3{(n+1)^2\bi{2n+2}{n+1}}-\f1{(n+1)^2}.
\endalign$$
Therefore (2.7) follows by induction. \qed
\medskip

The following lemma is essentially known, but we will include a simple proof.

\proclaim{Lemma 2.4} For any prime $p>3$ we have
$$\sum_{k=1}^{(p-1)/2}\f1{k^2}\eq0\ (\mo\ p)
\ \t{and}\ \sum_{k=1}^{(p-1)/2}\f{(-1)^k}{k^2}\eq (-1)^{(p-1)/2}\,2E_{p-3}\ (\mo\ p).$$
\endproclaim
\Proof. Since $\sum_{j=1}^{p-1}1/(2j)^2\eq\sum_{k=1}^{p-1}1/k^2\ (\mo\ p)$, we have the well-known congruence
$\sum_{k=1}^{p-1}1/k^2\eq0\ (\mo\ p)$. Thus
$$\sum_{k=1}^{(p-1)/2}\f1{k^2}\eq\f12\sum_{k=1}^{(p-1)/2}\l(\f1{k^2}+\f1{(p-k)^2}\r)=\f12\sum_{k=1}^{p-1}\f1{k^2}\eq0\ (\mo\ p).$$
By Lehmer [L, (20)],
$$\sum_{k=1}^{\lfloor p/4\rfloor}\f1{k^2}\eq(-1)^{(p-1)/2}\,4E_{p-3}\ (\mo\ p).$$
Therefore
$$\sum_{k=1}^{(p-1)/2}\f{1+(-1)^k}{k^2}=\sum_{j=1}^{\lfloor p/4\rfloor}\f2{(2j)^2}
\eq(-1)^{(p-1)/2}\,2E_{p-3}\ (\mo\ p)$$
and hence the second congruence in Lemma 2.4 also holds. \qed

\medskip
\noindent{\it Proof of (1.2)-(1.5)}. Note that (1.4) and (1.5) hold trivially when $p=3$.
Below we assume that $p=2n+1>3$.

With the help of Lemma 2.1, we have
$$\align&\sum_{k=n+1}^{p-1}\f{\bi{2k}k}k=\sum_{k=n+1}^{p-1}\f{k\bi{2k}k}{k^2}
\\\eq&\sum_{k=n+1}^{p-1}\f{2p}{k^2\bi{2(p-k)}{p-k}}=\sum_{j=1}^{n}\f{2p}{(p-j)^2\bi{2j}j}
\eq\sum_{k=1}^{n}\f{2p}{k^2\bi{2k}k}\ (\mo\ p^2).
\endalign$$
As $\sum_{k=0}^{p-1}\bi{2k}k/k\eq0\ (\mo\ p^2)$ by [ST1],
we obtain that
$$\sum_{k=1}^{n}\f{\bi{2k}k}k\eq-\sum_{k=n+1}^{p-1}\f{\bi{2k}k}k
\eq-2p\sum_{k=1}^{n}\f1{k^2\bi{2k}k}\ (\mo\ p^2).\tag2.8$$
In view of (2.5),
$$\sum_{k=1}^n\f{\bi{2k}k}k=\f{2n+1}3\bi{2n}n\sum_{k=1}^n\f1{k^2\bi nk^2}
\eq\f p3(-1)^n\sum_{k=1}^n\f1{k^2\bi nk^2}\ (\mo\ p^2).\tag2.9$$
Since
$$\bi{n+k}k(-1)^k=\bi{-n-1}k\eq\bi{p-n-1}k=\bi nk\ (\mo\ p)$$
for every $k=1,\ldots,n$, (2.6) yields that
$$\sum_{k=1}^n\f1{k^2\bi nk^2}\eq(-1)^{n-1}\(3\sum_{k=1}^n\f1{k^2\bi{2k}k}+2\sum_{k=1}^n\f{(-1)^k}{k^2}\)\ (\mo\ p).
\tag 2.10$$
Combining (2.8)--(2.10) we get
$$3\sum_{k=1}^n\f{\bi{2k}k}k+2p\sum_{k=1}^n\f{(-1)^k}{k^2}\eq-3p\sum_{k=1}^n\f1{k^2\bi{2k}k}
\eq\f 32\sum_{k=1}^n\f{\bi{2k}k}k\ (\mo\ p^2).\tag2.11$$
In view of Lemma 2.4,
$$\sum_{k=1}^n\f{(-1)^k}{k^2}\eq(-1)^n 2E_{p-3}\ (\mo\ p). $$
So, we have (1.2) and (1.3) by (2.11) and (2.8).

 By Lemma 2.1, (1.4) and (1.5) are equivalent.
Since $\sum_{k=1}^n1/k^2\eq0\ (\mo\ p)$ (by Lemma 2.4) and
 $$\bi{n+k}k\eq\bi{k-1/2}k=\f{\bi{2k}k}{4^k}\qquad \t{for every}\ k=1,\ldots,n,$$
 we obtain (1.4) from (2.7) and (1.3). \qed

\heading{3. Proof of Theorem 1.2}\endheading

\proclaim{Lemma 3.1} Let $p=2n+1$ be an odd prime. For $k=0,\ldots,n$ we have
$$\bi{n+k}{2k}-\f{p\bi nk}{4^{k+1}}(H_{n+k}-H_{n-k})\eq\f{\bi{2k}k}{(-16)^k}\ (\mo\ p^3)\tag3.1$$
and
$$\bi nk\bi{n+k}k(-1)^k\l(1-\f p4(H_{n+k}-H_{n-k})\r)\eq\f{\bi{2k}k^2}{16^k}\ (\mo\ p^4).\tag3.2$$
\endproclaim
\Proof. Both (3.1) and (3.2) hold trivially when $k=0$. Below we fix $k\in\{1,\ldots,n\}$.

As noted by the author's brother Z. H. Sun,
$$\align\bi{n+k}{2k}=&\f{\prod_{j=1}^k(p^2-(2j-1)^2)}{4^k\times(2k)!}
\\=&\f{\prod_{j=1}^k(-(2j-1)^2)}{4^k\times(2k)!}\prod_{j=1}^k\l(1-\f{p^2}{(2j-1)^2}\r)
\\\eq&\f{\bi{2k}k}{(-16)^k}\(1-\sum_{j=1}^k\f{p^2}{(2j-1)^2}\)\ (\mo\ p^4).
\endalign$$
Observe that
$$\align H_{n+k}-H_{n-k}=&H_n+\sum_{j=1}^k\f1{n+j}-H_n+\sum_{j=1}^k\f1{n+1-j}
\\=&\sum_{j=1}^k\f{2n+1}{(n+j)(n+1-j)}=\sum_{j=1}^k\f{p}{p^2/4-(j-1/2)^2}
\\\eq&-4\sum_{j=1}^k\f{p}{(2j-1)^2}\ (\mo\ p^3).
\endalign$$
Therefore $$\bi{n+k}{2k}\eq\f{\bi{2k}k}{(-16)^k}\l(1+\f p4(H_{n+k}-H_{n-k})\r)\ (\mo\ p^4).\tag3.3$$

Note that $p(H_{n+k}-H_{n-k})\eq0\ (\mo\ p^2)$ and
$$\bi nk\eq\bi{-1/2}k=\f{\bi{2k}k}{(-4)^k}\ (\mo\ p).$$
Thus (3.1) follows from (3.3) immediately.

In light of (3.3), we have
$$\align &\bi nk\bi{n+k}k(-1)^k=\bi{n+k}{2k}\bi{2k}k(-1)^k
\\\eq&\f{\bi{2k}k^2}{16^k}\l(1+\f p4(H_{n+k}-H_{n-k})\r)
\\\eq&\f{\bi{2k}k^2}{16^k}+\bi nk\bi{n+k}k(-1)^k\f p4(H_{n+k}-H_{n-k})\ (\mo\ p^4).
\endalign$$
(Recall that $p(H_{n+k}-H_{n-k})\eq0\ (\mo\ p^2)$.) So (3.2) also holds.
 \qed

\proclaim{Lemma 3.2} For any $n\in\Z^+$ we have
$$(-1)^n\sum_{k=0}^n\bi nk(-2)^{n-k}(H_{n+k}-H_{n-k})=\sum_{k=1}^n\f{(-1)^k}k-\f{H_n}2,\tag3.4$$
$$(-1)^n\sum_{k=0}^n\bi nk\bi{n+k}k(-1)^k(H_{n+k}-H_{n-k})=\f32\sum_{k=1}^n\f{\bi{2k}k}k\tag3.5$$
and
$$\aligned&(-1)^n\sum_{k=0}^n\bi nk\bi{n+k}k(-1)^k k(H_{n+k}-H_{n-k})
\\=&(2n+1)\l(1-\bi{2n}n\r)+\f 32 n(n+1)\sum_{k=1}^n\f{\bi{2k}k}k.
\endaligned\tag3.6$$
\endproclaim
\Proof. Via the software {\tt Sigma} we find the identities
$$\align\sum_{k=0}^n\bi nk(-2)^{n-k}H_{n+k}=&(-1)^n\f{H_n}2,\tag3.7
\\\sum_{k=0}^n\bi nk(-2)^{n-k}H_{n-k}=&(-1)^nH_n-(-1)^n\sum_{k=1}^n\f{(-1)^k}k,\tag3.8
\\\sum_{k=0}^n\bi nk\bi{n+k}k(-1)^kH_{n-k}=&2(-1)^nH_n-\f32(-1)^n\sum_{k=1}^n\f{\bi{2k}k}k.\tag3.9
\endalign$$
These identities can be easily proved by the WZ method (see, e.g., [PWZ]).
(The reader may consult [OS] to see how to produce such identities.)
Also,  it is known that (cf. [OS] and [Pr])
$$\sum_{k=0}^n\bi nk\bi{n+k}k(-1)^kH_{n+k}=(-1)^n2H_n.$$
By [OS, (36) and (37)], we have
$$(-1)^n\sum_{k=0}^n\bi nk\bi{n+k}k(-1)^k kH_{n+k}=2n(n+1)H_n-n^2$$
and
$$\align&(-1)^n\sum_{k=0}^n\bi nk\bi{n+k}k(-1)^k kH_{n-k}
\\=&(2n+1)\bi{2n}n-(n+1)^2+2n(n+1)H_n-\f32n(n+1)\sum_{k=1}^n\f{\bi{2k}k}k.
\endalign$$

In view of the above six identities we immediately obtain the desired (3.4)-(3.6). \qed

\Remark\ 3.1. S. Ahlgren and Ono [AO] employed the identity
$$\sum_{k=1}^n\bi nk^2\bi{n+k}k^2(1+2k(H_{n+k}+H_{n-k})-4kH_k)=0$$
to prove a super congruence conjectured by F. Beukers [Be].
\medskip

\noindent{\it Proof of Theorem 1.2}.
Set $n=(p-1)/2$.
In light of (3.1) and (3.4), we have
$$\align&\sum_{k=0}^n\f{\bi{2k}k}{8^k}-\sum_{k=0}^n(-2)^k\bi{n+k}{2k}
\\\eq&-\f p4\sum_{k=0}^n\f{\bi nk}{(-2)^k}(H_{n+k}-H_{n-k})
\\=&\f{p}{4\times2^n}\(\f{H_n}2-\sum_{k=1}^n\f{(-1)^k}k\)\ (\mo\ p^3).
\endalign$$
By a known identity (cf. (1.62) of [Go, p.8]),
$$\sum_{k=0}^n(-1)^k\bi{2n-k}k(2\cos x)^{2(n-k)}=\f{\sin((2n+1)x)}{\sin x}$$
and hence
$$\align\sum_{k=0}^n(-2)^k\bi{n+k}{2k}=&\sum_{k=0}^n(-1)^{n-k}\bi{2n-k}{k}\l(2\cos\f{\pi}4\r)^{2(n-k)}
\\=&(-1)^n\f{\sin((2n+1)\pi/4)}{\sin(\pi/4)}
=\l(\f2{2n+1}\r)=\l(\f2p\r).
\endalign$$
In view of (2.2) and (2.3), we also have
$$\align\f{H_n}2-\sum_{k=1}^n\f{(-1)^k}k=&\f32H_n-\sum_{k=1}^n\f{1+(-1)^k}k
\\=&\f32H_{(p-1)/2}-H_{\lfloor p/4\rfloor}
\\\eq&(-1)^{(p-1)/2}pE_{p-3}\ (\mo\ p^2).
\endalign$$
Therefore
$$\sum_{k=0}^n\f{\bi{2k}k}{8^k}-\l(\f2p\r)\eq\f p{4\times2^n}\l(\f{-1}p\r)pE_{p-3}
\eq\f{p^2}4\l(\f{-2}p\r)E_{p-3}\ (\mo\ p^3).$$
This proves (1.6).

By (3.2) and (3.5),
$$\sum_{k=0}^n\f{ \bi{2k}k^2}{16^k}-\sum_{k=0}^n\bi nk\bi{n+k}k(-1)^k
\eq-\f p4(-1)^n\f 32\sum_{k=1}^n\f{\bi{2k}k}k\ (\mo\ p^4).$$
With the help of the Chu-Vandermonde identity (cf. [GKP, p.\,169]),
$$\sum_{k=0}^n\bi nk\bi{n+k}k(-1)^k=\sum_{k=1}^n\bi n{n-k}\bi{-n-1}k=\bi{-1}n=(-1)^n.$$
Thus
$$\sum_{k=0}^n\f{\bi{2k}k^2}{16^k}-(-1)^n\eq-(-1)^n\f 38p\sum_{k=1}^n\f{\bi{2k}k}k\ (\mo\ p^4)$$
which gives (1.11).

By the Chu-Vandermonde identity we also have
$$\align&\sum_{k=0}^n k\bi nk\bi{n+k}k(-1)^k
\\=&n\sum_{k=1}^n\bi{n-1}{k-1}\bi{-n-1}k
=n\sum_{k=0}^n\bi{n-1}{n-k}\bi{-n-1}k
\\=&n\bi{-2}n=(-1)^n n(n+1).
\endalign$$
Recall Morley's congruence (cf. [Mo] and [P])
$$\bi{2n}n\eq(-1)^n4^{p-1}\ (\mo\ p^3).$$
Note also that $n(n+1)=(p^2-1)/4$ and
$$\sum_{k=1}^n\f{\bi{2k}k}k\eq\sum_{k=1}^{p-1}\f{\bi{2k}k}k\eq0\ (\mo\ p).$$
Therefore, by (3.2) and (3.6) we have
$$\align&\sum_{k=0}^n\f{k\bi{2k}k^2}{16^k}-(-1)^n\f{p^2-1}4
\\\eq&-\f p4\(p((-1)^n-4^{p-1})+(-1)^n\f 32\cdot\f{p^2-1}4\sum_{k=1}^n\f{\bi{2k}k}k\)
\\\eq&\f {p^2}4(4^{p-1}-(-1)^n)+\f 3{32}\l(\f{-1}p\r)p\sum_{k=1}^n\f{\bi{2k}k}k\ (\mo\ p^4)
\endalign$$
and hence
$$\align&\sum_{k=0}^n\f{k\bi{2k}k^2}{16^k}+\f{(-1)^n}4
\\\eq&\f {p^2}4(1+(2^{p-1}-1))^2+\f 3{32}\l(\f{-1}p\r)p\sum_{k=1}^n\f{\bi{2k}k}k
\\\eq&\f{p^2}4(2^{p}-1)+\f 3{32}\l(\f{-1}p\r)p\sum_{k=1}^n\f{\bi{2k}k}k\ (\mo\ p^4).
\endalign$$
This proves (1.12).

In light of (1.2), clearly (1.7) and (1.8) follow from (1.11) and (1.12) respectively.

Now we prove (1.9). By Lemma 2.1,
$$\align \f1{p^2}\sum_{p/2<k<p}\f{\bi{2k}k^2}{16^k}\eq&\sum_{p/2<k<p}\f{4/\bi{2(p-k)}{p-k}^2}{k^216^k}
=\sum_{k=1}^n\f 4{(p-k)^2\bi{2k}k^216^{p-k}}
\\\eq&\f14\sum_{k=1}^n\f{16^k}{k^2\bi{2k}k^2}\eq\f14\sum_{k=1}^n\f1{k^2\bi nk^2}\ (\mo\ p).
\endalign$$
This, together with (2.9) and (1.2), yields (1.9).

By Tauraso's identity mentioned in Remark 1.3,
$$\sum_{p/2<k<p}\f{4k+1}{16^k}\bi{2k}k^2=\f{p^2}{16^{p-1}}\bi{2p-1}{p-1}^2
-\f{p^2}{4^{p-1}}\bi{p-1}{(p-1)/2}^2\eq0\ (\mo\ p^3).$$
So (1.10) follows from (1.9).

\smallskip

The proof of Theorem 1.2 is now complete.  \qed

\medskip

\heading{4. Proof of Theorem 1.3}\endheading

 As we mentioned in the paragraph after Remark 1.4, based on the author's conjecture that $t_n\in\Z$
 for all $n\in\Z^+$, Kasper Andersen obtained the following lemma.

\proclaim{Lemma 4.1 {\rm (Kasper Andersen)}} For any $n\in\Z^+$ the number
$$t_n:=\f1{4n\bi{2n}n}\sum_{k=0}^{n-1}(21k+8)\bi{2k}k^3$$
is indeed an integer as conjectured by Z. W. Sun; in fact,
$$t_n=\sum_{k=0}^{n-1}\bi{n+k-1}k^2.$$
\endproclaim

We also need the following result.

\proclaim{Lemma 4.2} Let $p>3$ be a prime. Then
$$\sum_{k=1}^{p-1}\f{1+2pH_{k-1}}{k^2}\eq \f83pB_{p-3}\ (\mo\ p^2).$$
\endproclaim
\Proof. Note that
$$\sum_{k=1}^{p-1}\f1{(-2k)^3}\eq\sum_{j=1}^{p-1}\f1{j^3}\ (\mo\ p)$$
and hence $\sum_{k=1}^{p-1}1/k^3\eq0\ (\mo\ p)$ since $1-(-2)^3\not\eq0\ (\mo\ p)$.
By [ST1, (5.3)], we have
$$\sum_{k=1}^{p-1}\f{1+2pH_{k}}{k^2}\eq \f83pB_{p-3}\ (\mo\ p^2).$$
As $H_k=H_{k-1}+1/k$ for $k\in\Z^+$, the desired result follows. \qed

\medskip
\noindent{\it Proof of Theorem 1.3}. In light of Lemma 4.1,
$$\f1{p^a}\sum_{k=0}^{p^a-1}(21k+8)\bi{2k}k^3=4\bi{2p^a}{p^a}\sum_{k=0}^{p^a-1}\bi{p^a+k-1}k^2.$$
So we turn to determining $\f12\bi{2p^a}{p^a}$ and $\sum_{k=1}^{p^a-1}\bi{p^a+k-1}k^2$ modulo $p^4$.

By a result of Glaisher [G1,\,G2], if $p>3$ then
$$\f12\bi{2p}{p}=\bi{2p-1}{p-1}\eq1-\f23 p^3B_{p-3}\ (\mo\ p^4).$$
In view of [SD, Lemma 3.2],
$$\f12\bi{2p^{i+1}}{p^{i+1}}\eq\f12\bi{2p^{i}}{p^{i}}\ (\mo\ p^{2i+2})\quad\t{for every}\ i=1,2,3,\ldots.$$
Thus
$$\f12\bi{2p^a}{p^a}\eq\f12\bi{2p}{p}\eq\cases p^2+(-1)^{p-1}\ (\mo\ p^4)&\t{if}\ p\in\{2,3\},\\
1-\f23p^3B_{p-3}\ (\mo\ p^4)&\t{if}\ p>3.\endcases$$

Observe that
$$\align &\sum_{k=1}^{p^a-1}\bi{p^a+k-1}k^2=\sum_{k=1}^{p^a-1}\(\f{p^a}k\prod_{0<j<k}\l(1+\f{p^a}j\r)\)^2
\\\eq&\sum_{k=1}^{p^a-1}\l(\f{p^a}k\r)^2\prod_{0<j<k}\l(1+2\f{p^a}j\r)
\eq\sum_{k=1}^{p-1}\l(\f{p^a}{p^{a-1}k}\r)^2\prod_{0<j<k}\l(1+2\f{p^a}{p^{a-1}j}\r)
\\\eq&p^2\sum_{k=1}^{p-1}\f{1+2pH_{k-1}}{k^2}
\eq\cases p^2/(p-1)\ (\mo\ p^4)&\t{if}\ p\in\{2,3\},\\
 8p^3B_{p-3}/3\ (\mo\ p^4)&\t{if}\ p>3.\endcases
\endalign$$
(In the last step we apply Lemma 4.2.)

Combining the above we see that
$$\align&\f1{p^a}\sum_{k=0}^{p^a-1}(21k+8)\bi{2k}k^3
\\\eq &\cases 8(p^2+(-1)^{p-1})(p^2/(p-1)+1)\ (\mo\ p^4)&\t{if}\ p\in\{2,3\},
\\8(1-\f23p^3B_{p-3})(1+\f 83p^3B_{p-3})\ (\mo\ p^4)&\t{if}\ p>3,\endcases
\\\eq& 8(1+2p^3B_{p-3})\ \ (\mo\ p^4).
\endalign$$
This proves (1.15). \qed

\heading{5. More conjectures}\endheading

In 1914 S. Ramanujan [R] found the following curious identities
(see [BB], [B, pp.\,353-354] and [BBC] for more such series):
$$\sum_{k=0}^\infty(6k+1)\f{\bi{2k}k^3}{256^k}=\f4{\pi},
\  \ \sum_{k=0}^\infty(6k+1)\f{\bi{2k}k^3}{(-512)^k}=\f{2\sqrt2}{\pi},$$
and
$$\sum_{k=0}^\infty(42k+5)\f{\bi{2k}k^3}{4096^k}=\f{16}{\pi}.$$
(They are usually stated in terms of Gaussian hypergeometric series.)
For an odd prime $p$, L. van Hamme [vH] conjectured
the following $p$-adic analogues of the above three identities of Ramanujan:
$$\align\sum_{k=0}^{(p-1)/2}(6k+1)\f{\bi{2k}k^3}{256^k}\eq&p\l(\f{-1}p\r)\ (\mo\ p^4)\ \ \t{if}\ p>3,
\\\sum_{k=0}^{(p-1)/2}(6k+1)\f{\bi{2k}k^3}{(-512)^k}\eq&p\l(\f{-2}p\r)\ (\mo\ p^3),
\\ \sum_{k=0}^{(p-1)/2}(42k+5)\f{\bi{2k}k^3}{4096^k}\eq&5p\l(\f{-1}p\r)\ (\mo\ p^4).
\endalign$$
The first of these was recently shown by L. Long [Lo];
the second and the third remain open.

Motivated by Theorem 1.3 and Lemma 4.1, we propose the following conjecture.

\proclaim{Conjecture 5.1} {\rm (i)} For each $n=2,3,\ldots$ we have
$$\align &2n\bi{2n}n\ \bigg|\ \sum_{k=0}^{n-1}(3k+1)\bi{2k}k^3(-8)^{n-1-k},
\\&2n\bi{2n}n\ \bigg|\ \sum_{k=0}^{n-1}(3k+1)\bi{2k}k^316^{n-1-k},
\\&2n\bi{2n}n\ \bigg|\ \sum_{k=0}^{n-1}(6k+1)\bi{2k}k^3256^{n-1-k},
\\&2n\bi{2n}n\ \bigg|\ \sum_{k=0}^{n-1}(6k+1)\bi{2k}k^3(-512)^{n-1-k},
\\&2n\bi{2n}n\ \bigg|\ \sum_{k=0}^{n-1}(42k+5)\bi{2k}k^34096^{n-1-k}.
\endalign$$

{\rm (ii)} Let $p>3$ be a prime. Then
$$\align\sum_{k=0}^{p-1}\f{3k+1}{(-8)^k}\bi{2k}k^3\eq& p\l(\f{-1}p\r)+p^3E_{p-3}\ (\mo\ p^4),
\\\sum_{k=0}^{(p-1)/2}\f{3k+1}{16^k}\bi{2k}k^3\eq& p+2\l(\f{-1}p\r)p^3E_{p-3}\ (\mo\ p^4),
\\\sum_{k=0}^{p-1}\f{6k+1}{256^k}\bi{2k}k^3\eq& p\l(\f{-1}p\r)-p^3E_{p-3}\ (\mo\ p^4),
\\\sum_{k=0}^{(p-1)/2}\f{6k+1}{(-512)^k}\bi{2k}k^3\eq&p\l(\f{-2}p\r)+\f{p^3}4\l(\f2p\r)E_{p-3}\ (\mo\ p^4),
\\\sum_{k=0}^{p-1}\f{42k+5}{4096^k}\bi{2k}k^3\eq& 5p\l(\f{-1}p\r)-p^3E_{p-3}\ (\mo\ p^4),
\endalign$$
and
$$\sum_{k=0}^{(p-1)/2}\f{3k+1}{(-8)^k}\bi{2k}k^3\eq 4\l(\f{2}p\r)\sum_{k=0}^{p-1}\f{6k+1}{(-512)^k}\bi{2k}k^3-3p\l(\f{-1}p\r)\ (\mo\ p^4).$$
Also, for any $a\in\Z^+$ we have
$$\f1{p^a}\sum_{k=0}^{p^a-1}\f{3k+1}{16^k}\bi{2k}k^3\eq1+\f 76p^3B_{p-3}\ (\mo\ p^4)$$
and
$$\f1{p^a}\sum_{k=0}^{(p^a-1)/2}\f{42k+5}{4096^k}\bi{2k}k^3\eq\l(\f{-1}{p^a}\r)\l(5-\f 34pH_{p-1}\r)\ (\mo\ p^5).$$
\endproclaim

Each of Ishikawa [I], van Hamme [vH], Ahlgren [A] and
Mortenson [M05] confirmed the following
conjecture of Rodriguez-Villegas via certain advanced tools:
$$\sum_{k=0}^{p-1}\f{\bi{2k}k^3}{64^k}\eq a(p)\ (\mo\ p^2)\ \ \t{for any odd prime}\ p,$$
where the sequence $\{a(n)\}_{n\gs1}$ is defined by
$$\sum_{n=1}^\infty a(n)q^n=q\prod_{n=1}^\infty(1-q^{4n})^6\ \ \ (|q|<1)$$
and related to the Dedekind $\eta$-function in the theory of modular forms.
In 1892 F. Klein and R. Fricke proved that (cf. [SB, Theorem 14.2])
$$a(p)=\cases 4x^2-2p&\t{if}\  p=x^2+y^2\ \t{with}\ 2\nmid x\ \t{and}\ 2\mid y,
\\0&\t{if}\ p\eq3\ (\mo\ 4).\endcases$$

Let $p$ be an odd prime. Since $\bi{-1/2}k=\bi{2k}k/(-4)^k\ (k=0,1,2,\ldots)$,
for any integer $x\not\eq0\ (\mo\ p)$ we have
$$\align\sum_{k=0}^{p-1}\f{\bi{2k}k^3}{64^k}x^k&
\eq\sum_{k=0}^{(p-1)/2}\bi{(p-1)/2}k^3(-x)^k
\\&=\sum_{j=0}^{(p-1)/2}\bi{(p-1)/2}j^3(-x)^{(p-1)/2-j}
\\&\eq\sum_{k=0}^{p-1}\f{\bi{2k}k^3}{64^k}\l(\f{-x}p\r)x^{-k}\ (\mo\ p).
\endalign$$
Via computation we find that
$$\sum_{k=0}^{p-1}\f{\bi{2k}k^3}{64^k}\l(x^k-\l(\f{-x}p\r)x^{-k}\r)\eq0\ (\mo\ p^2)$$
for $x=1,4,-8,64$. (Note that the case $x=1$ is clear.) This leads us to propose the following conjecture.

\proclaim{Conjecture 5.2} Let $p$ be an odd prime.

{\rm (i)} If $p\eq1\pmod4$, then
$$\sum_{k=0}^{p-1}\f{\bi{2k}k^3}{(-8)^k}\eq\sum_{k=0}^{p-1}\f{\bi{2k}k^3}{64^k}
\eq\l(\f{2}p\r)\sum_{k=0}^{p-1}\f{\bi{2k}k^3}{(-512)^k}\pmod{p^3};$$
if $p\eq3\pmod4$ then
$$\sum_{k=0}^{p-1}\f{\bi{2k}k^3}{(-8)^k}\eq\sum_{k=0}^{p-1}\f{\bi{2k}k^3}{(-512)^k}\eq0\pmod{p^2}.$$

{\rm (ii)} If $p\eq1\pmod3$ and $p=x^2+3y^2$ with $x,y\in\Z$, then
$$\align\sum_{k=0}^{p-1}\f{\bi{2k}k^3}{16^k}\eq&\l(\f{-1}p\r)\sum_{k=0}^{p-1}\f{\bi{2k}k^3}{256^k}\pmod{p^3}
\\\eq&4x^2-2p\pmod{p^2};
\endalign$$
if $p\eq2\pmod3$, then
$$\sum_{k=0}^{p-1}\f{\bi{2k}k^3}{16^k}\eq\sum_{k=0}^{p-1}\f{\bi{2k}k^3}{256^k}\eq0\pmod{p^2}.$$

{\rm (iii) If $p\eq 1,2,4\pmod7\ (\t{i.e.,}\ (\f p7)=1)$, then
$$\sum_{k=0}^{p-1}\bi{2k}k^3\eq\l(\f{-1}p\r)\sum_{k=0}^{p-1}\f{\bi{2k}k^3}{4096^k}\ (\mo\ p^3);$$
if $p\eq3,5,6\pmod7\ (\t{i.e.,}\ (\f p7)=-1)$, then
$$\sum_{k=0}^{p-1}\f{\bi{2k}k^3}{4096^k}\eq0\pmod{p^2}.$$

{\rm (iv)} We have
$$\sum_{k=0}^{p-1}\f{\bi{2k}k^3}{(-64)^k}\eq\cases(\f{-1}p)(4x^2-2p)\pmod{p^2}&\t{if}\ (\f{-2}p)=1\ \&\ p=x^2+2y^2\ (x,y\in\Z),
\\0\pmod{p^2}&\t{if}\ (\f{-2}p)=-1,\ \t{i.e.},\ p\eq5,7\pmod8.\endcases$$
\endproclaim
\Remark\ 5.1. Let $p$ be an odd prime.
By the theory of binary quadratic forms (cf. Cox [C]), if $p\eq1\pmod3$ then there are unique $x,y\in\Z^+$
 such that $p=x^2+3y^2$; if $p\eq1,3\pmod8$ (i.e., $(\f{-2}p)=1$) then there are unique
$x,y\in\Z^+$ such that $p=x^2+2y^2$.

\proclaim{Conjecture 5.3} Let $p>3$ be a prime. Then
$$\sum_{k=0}^{p-1}\f{\bi{2k}k^2\bi{4k}{2k}}{81^k}\eq\cases\sum_{k=0}^{p-1}\bi{2k}k^3\pmod{p^3}&\t{if}\ (\f p7)=1,
\\0\pmod{p^2}&\t{if}\ (\f p7)=-1.\endcases$$
Also,
$$\f1{p^a}\sum_{k=0}^{(p^a-1)/2}\f{35k+8}{81^k}\bi{2k}k^2\bi{4k}{2k}\eq 8\times3^{p-1}\ (\mo\ p^2)$$
and
$$\f1{p^a}\sum_{k=0}^{p^a-1}\f{35k+8}{81^k}\bi{2k}k^2\bi{4k}{2k}\eq 8+\f{416}{27}p^3B_{p-3}\ (\mo\ p^4)$$
for all  $a\in\Z^+$.
Furthermore, for each $n=1,2,3,\ldots$ we have
$$\f1{4n(2n+1)\bi{2n}n}\sum_{k=0}^{n-1}(35k+8)\bi{2k}k^2\bi{4k}{2k} 81^{n-1-k}\in 3^{-\da(2n+1)}\Z,$$
where $\da(m)$ takes $1$ or $0$ according as $m$ is a power of $3$
or not.
\endproclaim

The author [S11a] made a conjecture on $\sum_{k=0}^{p-1}\bi{2k}k^2\bi{3k}k/64^k$ mod $p^2$ for any odd prime $p$.
Here we give a related conjecture.

\proclaim{Conjecture 5.4}
{\rm (i)} For any odd prime $p$ and positive integer $a$, we have
$$\f1{p^a}\sum_{k=0}^{p^a-1}\f{11k+3}{64^k}\bi{2k}k^2\bi{3k}k\eq3+\f 72p^3B_{p-3}\ (\mo\ p^4).$$
Moreover,
$$\f1{n(2n+1)\bi{2n}n}\sum_{k=0}^{n-1}(11k+3)\bi{2k}k^2\bi{3k}k64^{n-1-k}\in\Z$$
for all $n=2,3,\ldots$.

{\rm (ii)} If $p>3$ is a prime, then
$$p\sum_{k=1}^{(p-1)/2}\f{(11k-3)64^k}{k^3\bi{2k}k^2\bi{3k}k}\eq 32q_p(2)-\f{64}3 p^2B_{p-3}\ (\mo\ p^3).$$
\endproclaim

\proclaim{Conjecture 5.5} Let $p$ be an odd prime.
Then
 $$\align&\sum_{k=0}^{p-1}\f{\bi{2k}k^2\bi{3k}k}{8^k}
 \\\eq&\cases
 4x^2-2p\ (\mo\ p^2)&\t{if}\ (\f {-2}p)=1\ \&\ p=x^2+2y^2\ (x,y\in\Z),
 \\0\ (\mo\ p^2)&\t{if}\ (\f {-2}p)=-1.
 \endcases\endalign$$
 Also, for any $a\in\Z^+$ we have
$$\f1{p^a}\sum_{k=0}^{p^a-1}\f{10k+3}{8^k}\bi{2k}k^2\bi{3k}k\eq3+\f{49}8p^3B_{p-3}\ (\mo\ p^4).$$
Moreover, for each $n=2,3,\ldots$ we have
$$\f1{n(2n+1)\bi{2n}n}\sum_{k=0}^{n-1}(10k+3)\bi{2k}k^2\bi{3k}k8^{n-1-k}\in\Z.$$
\endproclaim

For $n\in\N$ the Bernoulli polynomial of degree $n$ is given by
$$B_n(x)=\sum_{k=0}^n\bi nkB_kx^{n-k}.$$

\proclaim{Conjecture 5.6} Let $p>3$ be a prime. Then
 $$\align&\sum_{k=0}^{p-1}\f{\bi{2k}k^2\bi{3k}k}{(-27)^k}
 \\\eq&\cases
 4x^2-2p\ (\mo\ p^2)&\t{if}\ p\eq1,4\ (\mo\ 15)\ \&\ p=x^2+15y^2\ (x,y\in\Z),
 \\2p-12x^2\ (\mo\ p^2)&\t{if}\ p\eq2,8\ (\mo\ 15)\ \&\ p=3x^2+5y^2\ (x,y\in\Z),
 \\0\ (\mo\ p^2)&\t{if}\ (\f p{15})=-1;
 \endcases\endalign$$
 $$\align&\sum_{k=0}^{p-1}\f{\bi{2k}{k}^2\bi{3k}{k}}{(-192)^k}
 \\\eq&\cases
 x^2-2p\ (\mo\ p^2)&\t{if}\ p\eq1\ (\mo\ 3)\ \&\ 4p=x^2+27y^2\ (x,y\in\Z),
 \\0\ (\mo\ p^2)&\t{if}\  p\eq2\ (\mo\ 3);
 \endcases\endalign$$
 $$\align&\sum_{k=0}^{p-1}\f{\bi{2k}k^2\bi{3k}k}{216^k}
 \eq\l(\f p3\r)\sum_{k=0}^{p-1}\f{\bi{2k}k^2\bi{4k}{2k}}{48^{2k}}\ \ \l(\mo\ p^{(5+(\f{-6}p))/2}\r)
 \\\eq&\cases
 4x^2-2p\ (\mo\ p^2)&\t{if}\ p\eq1,7\ (\mo\ 24)\ \&\ p=x^2+6y^2\ (x,y\in\Z),
 \\8x^2-2p\ (\mo\ p^2)&\t{if}\ p\eq5,11\ (\mo\ 24)\ \&\ p=2x^2+3y^2\ (x,y\in\Z),
 \\0\ (\mo\ p^2)&\t{if}\ (\f {-6}p)=-1,\ i.e.,\ p\eq 13,17,19,23\ (\mo\ 24);
 \endcases\endalign$$
$$\sum_{k=0}^{p-1}\f{\bi{2k}k^2\bi{4k}{2k}}{28^{4k}}\eq\cases 4x^2-2p\pmod{p^2}&\t{if}\ (\f{-2}p)=1\ \&\ p=x^2+2y^2\ (x,y\in\Z),
\\0\pmod{p^2}&\t{if}\ p\not=7\ \t{and}\ p\eq5,7\pmod8;
\endcases$$
$$\align&\sum_{k=0}^{p-1}\f{\bi{2k}k^2\bi{4k}{2k}}{(-2^{12}3)^k}
 \\\eq&\cases
 4x^2-2p\ (\mo\ p^2)&\t{if}\ 12\mid p-1,\ p=x^2+y^2,\ 3\nmid x\ \t{and}\ 3\mid y,
 \\-(\f{xy}3)4xy\ (\mo\ p^2)&\t{if}\ 12\mid p-5\ \t{and}\ p=x^2+y^2\  (x,y\in\Z),
 \\0\ (\mo\ p^2)&\t{if}\ p\eq3\ (\mo\ 4).
 \endcases\endalign$$
 Also, for any $a\in\Z^+$ we have
$$\align\f1{p^a}\sum_{k=0}^{p^a-1}\f{15k+4}{(-27)^k}\bi{2k}k^2\bi{3k}k\eq& 4\l(\f {p^a}3\r)+\l(\f{p^{a-1}}3\r)\f43p^2B_{p-2}\l(\f13\r)\ (\mo\ p^3),
\\\f1{p^a}\sum_{k=0}^{p^a-1}\f{5k+1}{(-192)^k}\bi{2k}{k}^2\bi{3k}{k}
\eq& \l(\f{p^a}3\r)+\l(\f{p^{a-1}}3\r)\f5{18}p^2B_{p-2}\l(\f13\r)\ (\mo\ p^3),
\\\f1{p^a}\sum_{k=0}^{p^a-1}\f{6k+1}{216^{k}}\bi{2k}k^2\bi{3k}k
 \eq& \l(\f {p^a}3\r)-\l(\f{p^{a-1}}3\r)\f5{12}p^2B_{p-2}\l(\f13\r)\ (\mo\ p^3),
\\\f1{p^a}\sum_{k=0}^{p^a-1}\f{8k+1}{48^{2k}}\bi{2k}k^2\bi{4k}{2k}\eq& \l(\f {p^a}3\r)-\l(\f{p^{a-1}}3\r)\f5{24}p^2B_{p-2}\l(\f13\r)\ (\mo\ p^3),
\\\f1{p^a}\sum_{k=0}^{p^a-1}\f{40k+3}{28^{4k}}\bi{2k}k^2\bi{4k}{2k}\eq& 3\l(\f {p^a}3\r)-\l(\f{p^{a-1}}3\r)\f{5p^2}{392}B_{p-2}\l(\f13\r)\, (\mo\ p^3)
\ \t{if}\ p\not=7,
\\\f1{p^a}\sum_{k=0}^{p^a-1}\f{28k+3}{(-2^{12}3)^k}\bi{2k}k^2\bi{4k}{2k}
\eq& 3\l(\f {p^a}3\r)+\l(\f{p^{a-1}}3\r)\f5{24}p^2B_{p-2}\l(\f13\r)\ (\mo\ p^3).
\endalign$$
and
$$\f 1{p^a}\sum_{p^a/2<k<p^a}\f{8k+1}{48^{2k}}\bi{2k}k^2\bi{4k}{2k}\eq0\pmod{p^2}.$$
 \endproclaim
 \Remark\ 5.2. (i) In view of Conjecture 5.6, we also have conjectures such as
$$\f1{2n(2n+1)\bi{2n}n}\sum_{k=0}^{n-1}(15k+4)\bi{2k}k^2\bi{3k}k(-27)^{n-1-k}\in 3^{-\da(2n+1)}\Z$$
and
$$\f1{2n(2n+1)\bi{2n}n}\sum_{k=0}^{n-1}(40k+3)\bi{2k}k^2\bi{4k}{2k}28^{4(n-1-k)}\in\Z,$$
where $n$ is any integer greater than one.
In addition, we guess that for any prime $p>3$ and $a\in\Z^+$ we have
 $$\f1{p^a}\sum_{k=0}^{p^a-1}\f{5k+1}{(-144)^k}\bi{2k}k^2\bi{4k}{2k}
\eq\l(\f{p^a}3\r)+\l(\f{p^{a-1}}3\r)\f 5{12}p^2B_{p-2}\l(\f13\r)\ (\mo\ p^3).$$

 (ii) The following Ramanujan-type series are closely related to some congruences in Conj. 5.6.
$$\gather\sum_{k=0}^\infty\f{5k+1}{(-192)^k}\bi{2k}k^2\bi{3k}k=\f{4\sqrt3}{\pi},
\ \sum_{k=0}^\infty\f{6k+1}{216^k}\bi{2k}k^2\bi{3k}k=\f{3\sqrt3}{\pi},
\\\sum_{k=0}^\infty\f{8k+1}{48^{2k}}\bi{2k}k^2\bi{4k}{2k}=\f{2\sqrt3}{\pi},
\ \ \sum_{k=0}^\infty\f{40k+3}{28^{4k}}\bi{2k}k^2\bi{4k}{2k}=\f{49}{3\sqrt3\,\pi},
\endgather$$ and
$$\sum_{k=0}^\infty\f{28k+3}{(-2^{12}3)^k}\bi{2k}k^2\bi{4k}{2k}=\f{16}{\sqrt3\,\pi}.$$
\medskip

 For the sake of brevity, below we will omit remarks like Remarks 5.1 and 5.2.

 For $n\in\N$ the Euler polynomial of degree $n$ is given by
$$E_n(x)=\sum_{k=0}^n\bi{n}{k}\frac{E_k}{2^k}\left(x-\frac{1}{2}\right)^{n-k}.$$

\proclaim{Conjecture 5.7} Let $p$ be an odd prime.
Then
 $$\align&\sum_{k=0}^{p-1}\f{\bi{2k}k^2\bi{4k}{2k}}{(-2^{10})^k}
 \\\eq&\cases
 4x^2-2p\ (\mo\ p^2)&\t{if}\ p\eq1,9\ (\mo\ 20)\ \&\ p=x^2+5y^2\ (x,y\in\Z),
 \\2p-2x^2\ (\mo\ p^2)&\t{if}\ p\eq3,7\ (\mo\ 20)\ \&\ 2p=x^2+5y^2\ (x,y\in\Z),
 \\0\ (\mo\ p^2)&\t{if}\ (\f {-5}p)=-1,\ i.e.,\ p\eq 11,13,17,19\ (\mo\ 20).
 \endcases\endalign$$
 and
 $$\f1{p^a}\sum_{k=0}^{p^a-1}\f{20k+3}{(-2^{10})^k}\bi{2k}k^2\bi{4k}{2k}\eq3\l(\f{-1}{p^a}\r)+3\l(\f{-1}{p^{a-1}}\r)p^2E_{p-3}\ (\mo\ p^3)
 \ \ \t{for all}\ a\in\Z^+.$$
 Provided $p>3$, we have
 $$\align&\sum_{k=0}^{p-1}\f{\bi{2k}k^2\bi{4k}{2k}}{12^{4k}}
 \\\eq&\cases
 4x^2-2p\ (\mo\ p^2)&\t{if}\ p\eq1,9,11,19\ (\mo\ 40)\ \&\ p=x^2+10y^2\ (x,y\in\Z),
 \\2p-8x^2\ (\mo\ p^2)&\t{if}\ p\eq7,13,23,37\ (\mo\ 40)\ \&\ p=2x^2+5y^2\ (x,y\in\Z),
 \\0\ (\mo\ p^2)&\t{if}\ (\f {-10}p)=-1,\ i.e.,\ p\eq 3,17,21,27,29,31,33,39\ (\mo\ 40),
 \endcases\endalign$$
 and
 $$\f1{p^a}\sum_{k=0}^{p^a-1}\f{10k+1}{12^{4k}}\bi{2k}k^2\bi{4k}{2k}
\eq \l(\f{-2}{p^a}\r)-\l(\f{-2}{p^{a-1}}\r)\f{p^2}{48}E_{p-3}\l(\f14\r)\ (\mo\ p^3)$$
for all $a\in\Z^+$.
When $p>5$, we have
$$\align &\sum_{k=0}^{p-1}\f{\bi{2k}k^2\bi{4k}{2k}}{(-2^{14}3^4 5)^k}
\\\eq&\cases4x^2-2p\ (\mo\ p^2)&\t{if}\ p\eq1,9\ (\mo\ 20),\ p=x^2+y^2,\ 5\nmid x\ \t{and}\ 5\mid y,
\\4xy\ (\mo\ p^2)&\t{if}\ p\eq13,17\ (\mo\ 20),\ p=x^2+y^2\ \t{and}\ 5\mid x+y,
\\0\ (\mo\ p^2)&\t{if}\ p\eq3\ (\mo\ 4).\endcases
\endalign$$
\endproclaim

\proclaim{Conjecture 5.8} Let $p$ be an odd prime. Then
$$\align&\sum_{k=0}^{p-1}\f{\bi{6k}{3k}\bi{3k}{k}\bi{2k}k}{(-2^{15})^k}
 \\\eq&\cases
 (\f{-2}p)(x^2-2p)\ (\mo\ p^2)&\t{if}\ (\f p{11})=1\ \&\ 4p=x^2+11y^2\ (x,y\in\Z),
 \\0\ (\mo\ p^2)&\t{if}\ (\f p{11})=-1,\ i.e.,\ p\eq2,6,7,8,10\ (\mo\ 11).
 \endcases\endalign$$
Also, for any $a\in\Z^+$ we have
$$\align &\f1{p^a}\sum_{k=0}^{p^a-1}\f{154k+15}{(-2^{15})^k}\bi{6k}{3k}\bi{3k}{k}\bi{2k}k
\\\eq& 15\l(\f {-2}{p^a}\r)
+\l(\f{-2}{p^{a-1}}\r)\f{15}{16}p^2E_{p-3}\l(\f14\r)\ (\mo\ p^3).
\endalign$$
Moreover, for each $n=2,3\ldots$ we have
$$\f{1}{2n(2n+1)\bi{2n}n}\sum_{k=0}^{n-1}(154k+15)\bi{6k}{3k}\bi{3k}{k}\bi{2k}k(-2^{15})^{n-1-k}\in\Z.$$
 \endproclaim

\proclaim{Conjecture 5.9} Let $p$ be an odd prime and let $a\in\Z^+$.
If $p\eq1\pmod3$, then
$$\sum_{k=0}^{p^a-1}\f{9k+2}{108^k}\bi{2k}k^2\bi{3k}k\eq0\pmod{p^{2a}}.$$
If $p\eq1,3\ (\mo\ 8)$, then
$$\sum_{k=0}^{p^a-1}\f{16k+3}{256^k}\bi{2k}k^2\bi{4k}{2k}\eq0\ (\mo\ p^{2a+\da_{p,3}}).$$
If $p\eq1\ (\mo\ 4)$, then
$$\sum_{k=0}^{p^a-1}\f{4k+1}{64^k}\bi{2k}k^3\eq0\ (\mo\ p^{2a})$$
and
$$\sum_{k=0}^{p^a-1}\f{36k+5}{12^{3k}}\bi{6k}{3k}\bi{3k}{k}\bi{2k}k\eq 0\ (\mo\ p^{2a+\da_{p,5}}).$$
\endproclaim
\Remark\ 5.3. The reader may consult conjectures of Rodriguez-Villegas [RV] (see also [M05])
on $\sum_{k=0}^{p-1}\bi{2k}k^2\bi{3k}k/108^k$, $\sum_{k=0}^{p-1}\bi{2k}k^2\bi{4k}{2k}/256^k$
and $\sum_{k=0}^{p-1}\bi{6k}{3k}\bi{3k}k\bi{2k}k/12^{3k}$ modulo $p^2$, where $p>3$ is a prime.
\medskip

We will give more conjectures similar to Conjectures 5.1-5.9 in the forthcoming survey [S11b].
Now we turn to sums of products of two binomial coefficients.

\proclaim{Conjecture 5.10} {\rm (i)} For any prime $p\eq1\pmod3$ and positive integer $a$, we have
$$\sum_{k=0}^{p^a-1}\f{k\bi{2k}k\bi{3k}k}{54^k}\eq0\pmod{p^{a+1}}.$$

{\rm (ii)} Let $p\eq1,3\pmod8$ be a prime and write $p=x^2+2y^2$ with $x\eq1\ (\mo\ 4)$. Then
 $$\sum_{k=0}^{p-1}\f{\bi{2k}k\bi{4k}{2k}}{128^k}\eq (-1)^{\lfloor(p+5)/8\rfloor}\l(2x-\f p{2x}\r)\pmod{p^2}.$$
 Also, for any $a\in\Z^+$ we have
 $$\sum_{k=0}^{p^a-1}\f{k\bi{2k}k\bi{4k}{2k}}{128^k}\eq0\pmod{p^{a+1+\da_{p,3}}}.$$

 {\rm (iii)} Let $p\eq1\ (\mo\ 4)$ be a prime and write
$p=x^2+y^2$ with $x\eq1\ (\mo\ 4)$ and $y\eq0\ (\mo\ 2)$. Then
$$\sum_{k=0}^{p-1}\f{\bi{6k}{3k}\bi{3k}k}{864^k}\eq\cases(-1)^{\lfloor x/6\rfloor}(2x-p/(2x))\ (\mo\ p^2)&\t{if}\ p\eq1\ (\mo\ 12),
\\(\f{xy}3)(2y-p/(2y))\ (\mo\ p^2)&\t{if}\ p\eq5\ (\mo\ 12).
\endcases$$
Also, for any $a\in\Z^+$ we have
$$\sum_{k=0}^{p^a-1}\f{k\bi{6k}{3k}\bi{3k}k}{864^k}\eq0\ (\mo\ p^{a+1})
\ \ \t{and}\ \ \f1{5^{a+2}}\sum_{k=0}^{5^a-1}\f{k\bi{6k}{3k}\bi{3k}k}{864^k}\eq3\ (\mo\ 5).$$
\endproclaim

\proclaim{Conjecture 5.11} Let $p>3$ be a prime.

{\rm (i)} If $p\eq7\ (\mo\ 12)$ and $p=x^2+3y^2$ with $y\eq1\ (\mo\ 4)$, then
$$\sum_{k=0}^{p-1}\l(\f k3\r)\f{\bi{2k}k^2}{(-16)^k}\eq (-1)^{(p-3)/4}\l(4y-\f p{3y}\r)\ (\mo\ p^2)$$
and we can determine $y$ mod $p^2$ via the congruence
$$\sum_{k=0}^{p-1}\l(\f k3\r)\f{k\bi{2k}k^2}{(-16)^k}\eq (-1)^{(p+1)/4}y\ (\mo\ p^2).$$

{\rm (ii)} If $p\eq1\ (\mo\ 12)$, then
$$\sum_{k=0}^{p-1}\bi{p-1}k\l(\f k3\r)\f{\bi{2k}k^2}{16^k}\eq0\ (\mo\ p^2).$$
If $p\eq11\ (\mo\ 12)$, then
$$\sum_{k=0}^{p-1}\l(\f k3\r)\f{\bi{2k}k^2}{(-16)^k}\eq0\ (\mo\ p).$$
\endproclaim
\Remark\ 5.4. The author could prove that $\sum_{k=0}^{p-1}(\f k3)\bi{2k}k^2/(-16)^k\eq0\pmod{p^2}$ for any prime $p\eq1\pmod4$.
\medskip

\proclaim{Conjecture 5.12} Let $p$ be an odd prime and let $a\in\Z^+$.

{\rm (i)} We have
$$\sum_{k=0}^{p^a-1}\f{\bi{2k}k\bi{4k}{2k}}{64^k}\eq\l(\f{-2}{p^a}\r)-\l(\f{-2}{p^{a-1}}\r)\f{3p^2}{16}E_{p-3}\l(\f14\r)\pmod{p^3}$$
and
$$\sum_{k=0}^{p^a-1}\f{\bi{2k}kC_{2k}}{64^k}\eq\l(\f{-1}{p^a}\r)-\l(\f{-1}{p^{a-1}}\r)3p^2E_{p-3}\pmod{p^3},$$
where $C_n$ stands for the Catalan number $\f1{n+1}\bi{2n}n=\bi{2n}n-\bi{2n}{n+1}$.

{\rm (ii)} Suppose $p>3$. Then
$$\sum_{k=0}^{p^a-1}\f{\bi{6k}{3k}\bi{3k}k}{432^k}\eq\l(\f{-1}{p^a}\r)-\l(\f{-1}{p^{a-1}}\r)\f{25}9p^2E_{p-3}\ (\mo\ p^3)$$
and
$$\sum_{k=0}^{p^a-1}\f{\bi{6k}{3k}C_k^{(2)}}{432^k}\eq\l(\f {p^a}3\r)\ (\mo\ p^2),$$
where $$C_k^{(2)}=\f{\bi{3k}k}{2k+1}=\bi{3k}k-2\bi{3k}{k-1}$$
is a second-order Catalan number.

{\rm (iii)} Assume $p>3$. Then
$$\align\sum_{k=1}^{p^a-1}\f{\bi{2k}{k+1}\bi{3k}{k+1}}{27^k}\eq&2\l(\f {p^a}3\r)-7\ (\mo\ p),
\\\sum_{k=1}^{p^a-1}\f{\bi{2k}{k-1}\bi{3k}{k-1}}{27^k}\eq&\l(\f {p^a}3\r)-p^a\ (\mo\ p^2),
\\\sum_{k=0}^{p^a-1}\f{\bi{2k}k\bi{3k}{k}}{27^k}\eq&\l(\f{p^a}3\r)-\l(\f{p^{a-1}}3\r)\f{p^2}3B_{p-2}\l(\f13\r)\pmod{p^3},
\\\sum_{k=0}^{p^a-1}\f{\bi{2k}kC_k^{(2)}}{27^k}\eq&\l(\f {p^a}3\r)-\l(\f{p^{a-1}}3\r)\f23p^2B_{p-2}\l(\f13\r)\ (\mo\ p^3).
\endalign$$
Furthermore,
$$\sum_{k=0}^{p^a-1}(4k+1)\f{\bi{2k}kC_k^{(2)}}{27^k}\eq\l(\f {p^a}3\r)\ (\mo\ p^4).$$
\endproclaim
\Remark\ 5.5. For a prime $p>3$, Rodriguez-Villegas' conjecture (cf. [RV]) on
$$\sum_{k=0}^{p-1}\f{\bi{2k}k\bi{3k}k}{27^k},\ \sum_{k=0}^{p-1}\f{\bi{2k}k\bi{4k}{2k}}{64^k},\ \sum_{k=0}^{p-1}\f{\bi{6k}{3k}\bi{3k}k}{432^k}$$
modulo $p^2$ were proved by Mortenson [M03b]. By Gosper's algorithm (cf. [PWZ]) we find that
$$\sum_{k=0}^{n}\f{9k+2}{27^k}\bi{2k}k\bi{3k}{k}=\f{(3n+1)(3n+2)}{27^n}\bi{2n}n\bi{3n}n$$
and
$$\sum_{k=0}^n\f{36k+5}{432^k}\bi{6k}{3k}\bi{3k}k=\f{(6n+1)(6n+5)}{432^n}\bi{6n}{3n}\bi{3n}n.$$
\medskip

\proclaim{Conjecture 5.13} Let $p>3$ be a prime. Then
$$\sum_{k=0}^{p-1}\f{\bi{2k}k\bi{3k}k}{24^k}\eq\l(\f p3\r)\sum_{k=0}^{p-1}\f{\bi{2k}k\bi{3k}k}{(-216)^k}
\eq\cases\bi{2(p-1)/3}{(p-1)/3}\ (\mo\ p^2)&\t{if}\ p\eq1\ (\mo\ 3),
\\p/\bi{2(p+1)/3}{(p+1)/3}\ (\mo\ p^2)&\t{if}\ p\eq 2\ (\mo\ 3).\endcases$$
Also,
$$\sum_{k=0}^{p-1}\f{\bi{3k}kC_k}{24^k}\eq\f19\l(\f p3\r)\sum_{k=0}^{p-1}\f{\bi{3k}kC_k}{(-216)^k}\eq\f12\bi{2(p-(\f p3))/3}{(p-(\f p3))/3}\ (\mo\ p).$$
When $p\eq1\ (\mo\ 3)$ and $4p=x^2+27y^2$ with $x,y\in\Z$ and $x\eq2\ (\mo\ 3)$, we have
$$x\eq\sum_{k=0}^{p-1}\f{k+2}{24^k}\bi{2k}k\bi{3k}k\eq\sum_{k=0}^{p-1}\f{9k+2}{(-216)^k}\bi{2k}k\bi{3k}k\pmod{p^2}.$$
\endproclaim

\proclaim{Conjecture 5.14} Let $p>3$ be a prime.

{\rm (i)} We always have
 $$\sum_{k=0}^{p-1}\f{\bi{2k}k\bi{4k}{2k+1}}{48^k}\eq0\ (\mo\
 p^2).$$
If $p\eq1\ (\mo\ 3)$ and $p=x^2+3y^2$ with $x\eq1\ (\mo\ 3)$, then
 $$\gather\sum_{k=0}^{p-1}\f{\bi{2k}k\bi{4k}{2k}}{48^k}\eq 2x-\f p{2x}\ (\mo\ p^2),
 \\\sum_{k=0}^{p-1}\f{k+1}{48^k}\bi{2k}k\bi{4k}{2k}\eq x\ (\mo\ p^2).
 \endgather$$
If $p\eq2\ (\mo\ 3)$, then
$$\sum_{k=0}^{p-1}\f{\bi{2k}k\bi{4k}{2k}}{48^k}\eq\f{3p}{2\bi{(p+1)/2}{(p+1)/6}}\pmod{p^2}.$$

{\rm (ii)} If
$(\f p7)=1$ and $p=x^2+7y^2$ with $(\f x7)=1$, then
 $$\gather\sum_{k=0}^{p-1}\f{\bi{2k}k\bi{4k}{2k}}{63^k}\eq \l(\f p3\r)\l(2x-\f p{2x}\r)\ (\mo\ p^2),
 \\\sum_{k=0}^{p-1}\f{k+8}{63^k}\bi{2k}k\bi{4k}{2k}\eq 8\l(\f p3\r)x\ (\mo\ p^2).
  \endgather$$
If $(\f p7)=-1$, then
$$\sum_{k=0}^{p-1}\f{\bi{2k}k\bi{4k}{2k}}{63^k}\eq\sum_{k=0}^{p-1}\f{\bi{2k}k\bi{4k}{2k}^2}{63^k}\eq0\ (\mo\ p).$$

{\rm (iii)} If $p\eq1\ (\mo\ 4)$ and
   $p=x^2+y^2$ with $x\eq1\ (\mo\ 4)$ and $y\eq0\ (\mo\ 2)$, then
 $$\gather\sum_{k=0}^{p-1}\f{\bi{2k}k\bi{4k}{2k}}{72^k}\eq \l(\f 6p\r)\l(2x-\f p{2x}\r)\ (\mo\ p^2),
 \\\sum_{k=0}^{p-1}\f{1-k}{72^k}\bi{2k}k\bi{4k}{2k}\eq\l(\f 6p\r)x\ (\mo\ p^2).
 \endgather$$
 If $p\eq3\ (\mo\ 4)$, then
 $$\sum_{k=0}^{p-1}\f{\bi{2k}k\bi{4k}{2k}}{72^k}\eq\l(\f 6p\r)\f{2p}{3\bi{(p+1)/2}{(p+1)/4}}\pmod{p^2}.$$
 \endproclaim
 \Remark\ 5.6. Let $p>3$ be a prime. Concerning part (i) the author could prove that
$$\sum_{k=0}^{p-1}\f{\bi{2k}k\bi{4k}{2k}}{(-192)^k}
\eq\l(\f{-2}p\r)\sum_{k=0}^{p-1}\f{\bi{2k}k\bi{4k}{2k}}{48^k}\ (\mo\ p^2)$$
and $$\sum_{k=0}^{p-1}\f{k\bi{2k}k\bi{4k}{2k}}{(-192)^k}
\eq\f14\l(\f{-2}p\r)\sum_{k=0}^{p-1}\f{k\bi{2k}k\bi{4k}{2k}}{48^k}\ (\mo\ p^2).$$
We have similar things related to parts (ii) and (iii) of Conj. 5.14.
\medskip

Finally, we mention that we also have some conjectural super congruences involving quadratic polynomials
and sums of products of more than three binomial coefficients.
For example, inspired by the identity
$$\sum_{k=1}^\infty\f{(-1)^k(205k^2-160k+32)}{k^5\bi{2k}k^5}=-2\zeta(3)$$
due to T. Amdeberhan and D. Zeilberger [AZ], on April 4, 2010 the author (cf. [S10a]) conjectured that
$$\sum_{k=0}^{(p-1)/2}(205k^2+160k+32)(-1)^k\bi{2k}k^5\eq 32p^2+\f{896}3p^5B_{p-3}\ (\mo\ p^6)$$
for any prime $p>3$. Also, (1.22) and (1.23) in Conjecture 1.4 were motivated by the first and the second congruences
in our following conjecture.

\proclaim{Conjecture 5.15} {\rm (i)} For any odd prime $p$, we have
$$\gather\sum_{n=0}^{p-1}\f{18n^2+7n+1}{(-128)^n}\bi{2n}n^2\sum_{k=0}^n\bi{-1/4}k^2\bi{-3/4}{n-k}^2
\eq p^2\l(\f 2p\r)\pmod{p^3},
\\\sum_{n=0}^{p-1}\f{40n^2+26n+5}{(-256)^n}\bi{2n}n^2\sum_{k=0}^n\bi nk^2\bi{2k}k\bi{2(n-k)}{n-k}
\eq5p^2\pmod{p^3},
\\\sum_{n=0}^{p-1}\f{12n^2+11n+3}{(-32)^n}\sum_{k=0}^n\bi nk^4\bi{2k}k\bi{2(n-k)}{n-k}\eq3p^2+\f 74p^5B_{p-3}\pmod{p^6}.
\endgather$$

{\rm (ii)} If $p>3$ is a prime, then
$$\sum_{n=0}^{p-1}\f{3n^2+n}{16^n}\sum_{k=0}^n\bi nk^2\bi{2k}k\bi{2(n-k)}{n-k}\eq-4p^4q_p(2)+6p^5q_p(2)^2\pmod{p^6}.$$
For any integer $m>1$, we have
$$a_m:=\f1{2m^3(m-1)}\sum_{n=0}^{m-1}(3n^2+n)16^{m-1-n}\sum_{k=0}^n\bi nk^2\bi{2k}k\bi{2(n-k)}{n-k}\in\Z;$$
moreover, $a_m$ is odd if and only if $m$ is a power of two.
\endproclaim

\medskip

\Ack. The author is indebted to Mr. Kasper Andersen for proving the author's conjecture
that $\f1n\sum_{k=0}^{n-1}(21k+8)\bi{2k}k^3$ is always an integer divisible by $4\bi{2n}n$.
He also thanks Prof. E. Mortenson and R. Tauraso for helpful comments on the initial version of this paper.

\enddocument